# Pursuit-Evasion in Graphs: Zombies, Lazy Zombies and a Survivor[§]

*Prosenjit Bose*[δ] *Jean-Lou De Carufel*[σ] *Thomas Shermer*[λ]

ABSTRACT. The game of cops and robber is a well-studied pursuit-evasion game played on graphs. We study a variant of this game known as the zombies and survivor game. In this variant, the single survivor plays the role of the robber and attempts to escape from the zombies that play the role of the cops. The zombies are restricted to always follow an edge of a shortest path towards the robber on their turn. Let $c(G)$ (respectively $z(G)$) be the smallest number of cops (respectively zombies) required to catch the robber (respectively the survivor) on a graph $G$. We show that there exist outerplanar graphs and 2-connected outerplanar graphs $G$ on $n$ vertices that require $\Omega(n)$ zombies to catch the survivor. Since 2 cops are sufficient to catch a robber on an outerplanar graph, we have that $\frac{z(G)}{c(G)} = \Omega(n)$. Previously, the best known lower bound was $\Omega(\log n)$ zombies for outerplanar graphs. We show that there exist simple polygons whose visibility graph $G$ is such that $\frac{z(G)}{c(G)} = \Omega(n)$. We also show that there exist maximum-degree-3 outerplanar graphs $G$ such that $\frac{z(G)}{c(G)} = \Omega\left(\frac{n}{\log(n)}\right)$.

Then, we establish that *lazy zombies* (which are zombies that are allowed to stay still on their turn) are more powerful than plain zombies and less powerful than cops. We prove that 2 lazy zombies always win (in less than $2n$ rounds) on connected outerplanar graphs. Then we show that $k$ lazy zombies win after $O(n^{2k})$ rounds on connected graphs with treedepth $k$. Finally, we highlight a few implications stemming from this upper bound such as $(k+1)\log n$ lazy zombies win on connected graphs with treewidth $k$, $O(\sqrt{n})$ lazy zombies are always sufficient to win on connected planar graphs, $O(\sqrt{gn})$ lazy zombies win on connected graphs with genus $g$ and $O(h\sqrt{hn})$ lazy zombies win on all connected graphs $G$ with any excluded $h$-vertex minor $H$. Our results on lazy zombies still hold when an adversary chooses the initial positions of the zombies.

## 1 Introduction

The game of cops and robber was first introduced by Quillot in his doctoral thesis [38] and then independently by Nowakowski and Winkler [36]. In this pursuit-evasion game, a set of cops move along the edges of a connected and undirected graph $G$ to catch a robber that is also moving along the edges of $G$. At the beginning, each cop chooses a starting vertex (multiple cops can occupy the same vertex). Then the robber chooses a starting vertex. From there, when it is the cops' turn to play, each cop decides either to stay still or to move to an adjacent vertex. When it is the robber's turn to play, the robber decides either to stay still or to move to an adjacent vertex. (In the classical version of the game, the cops and the robber are aware of each others' locations at all times.). If at least one cop reaches a vertex where the robber is standing, then the cops win. However, if the robber can escape for an infinite number of turns, then the robber wins. The *cop number* of a graph $G$, denoted as $c(G)$, is the minimum number of cops required to catch the robber on $G$. We refer the reader to the book by Bonato and Nowakowski for an excellent introduction to this topic [13].

[§]Research supported in part by NSERC.
[δ]Carleton University, jit@scs.carleton.ca
[σ]University of Ottawa, jdecaruf@uottawa.ca
[λ]Simon Fraser University, shermer@sfu.ca

Cops and robber has been studied on several classes of graphs. Quillot [38], and Nowakowski and Winkler [36] characterized graphs with cop number 1, which we call *cop-win graphs*. Aigner and Fromme showed that the cop number of planar graphs is 3 [1]. Clarke showed that the cop number of outerplanar graphs is 2 [16]. More generally, Schröder showed that the cop number of graphs with *genus g* is at most $\lfloor \frac{3}{2}g \rfloor + 3$ [40]. Given a graph $G$ with treewidth $k$, Joret et al. proved that the cop-number of $G$ is at most $\frac{1}{2}k + 1$ [27]. Andreae proved an upper bound on the cop number of graphs with excluded minor [4]. More specifically, let $H$ be a graph and $v$ be a vertex of $H$ such that $H - v$ has no isolated vertices. If $H$ is not a minor of $G$, then $c(G) < |E(H - v)|$. If we assume that $H$ is connected, it has a vertex $v$ such that $H - v$ has no isolated vertices, and we get $c(G) < \frac{1}{2}(h-1)(h-2)$, where $h$ is the number of vertices in $H$. Meyniel conjectured that in general, $c(G) = O(\sqrt{n})$.

Cops and robbers have also been studied on geometric graphs. Typically, a geometric graph is a graph whose vertices are points in the plane and whose edges are segments joining vertices that are often weighted by their length. Beveridge et al. showed that the cop number of unit disk graphs is between 3 and 9 [8]. Gavenciak et al. provide upper bounds on the cop numbers of different intersection graphs [23]. Lubiw et al. showed that the visibility graph of a simple polygon $P$ is cop-win [33]. They extend this result to a version of the game where players are allowed to move along any straight line segments inside $P$. In both versions, a direct consequence of their proof is that the cop catches the robber within $O(n_r)$ steps, where $n_r$ is the number of reflex vertices in $P$.

Multiple variants of the game of cops and robbers have been studied. For instance, the cop can capture the robber from an integer distance $k \geq 0$ [10, 11], a number of players can move with different speeds [21, 14, 2], or the robber can be invisible or even drunk [28, 29, 30]. The game of *active cops and robber*, where the robber must move at each of its turns, was first introduced by Aigner and Fromme [1], and later studied by Neufeld and Nowakowski [35]. Offner and Ojakian introduced a class of Cops and Robbers variants, where one has to specify how many cops must move on every cop turn, how many must remain in place, and how many may do either [37]. In this class of variants, one can define a version called *Lazy Cops and Robber*, where only one cop may move on each turn. This was studied by various authors [37, 5, 6, 22]. In the *Fully Active Cops and Robber*, no player is allowed to remain on a vertex [25].

Zombies and Survivor is a variant of the game of Cops and Robber. The deterministic version of Zombies and Survivor was first introduced by Fitzpatrick et al. [20]. In this game, each cop (thought of as a *zombie*) is restricted to move along an edge incident to its current position and belonging to a shortest path to the robber (thought of as the *survivor*). Moreover, the zombie is active, in the sense that it must move on its turn. If there is more than one shortest path between a given zombie and the survivor, the zombie can choose which path to follow. Since the zombies are active, the starting configuration of the zombies actually plays a role in determining whether a survivor is caught. For example, if several zombies are placed on the same vertex in a cycle with more than 4 vertices, then they will never catch a survivor. However, 2 zombies can be strategically placed in order to always catch a survivor on such a cycle. As such, we define two types of zombie number, one where the zombies are strategically placed and the other where an adversary determines the initial position of the zombies. The *zombie number* $z(G)$ of a graph $G$ is then defined as the minimum number of zombies required to catch the survivor on $G$, and the *universal zombie number* $u(G)$ is defined as the minimum number of zombies required to catch the survivor when the starting configuration of the zombies is determined by an adversary. The above example shows that $u(G) = \infty$ when $G$ is a cycle on more than 4 vertices. Since a cop has more power than a zombie, we have $c(G) \leq z(G) \leq u(G)$. From this observation, we get that zombie-win graphs are also cop-win graphs.

The randomized version of Zombies and Survivor was first introduced by Bonato et al. [12]. In this version, the zombies choose their initial starting locations uniformly at random. Then, on its turn, when a zombie has several edges that it can follow to get closer to the survivor, the zombie chooses one uniformly at random. In the randomized version, the *zombie number* of a graph is defined as the minimum number of zombies required such that the probability that they win is at least $\frac{1}{2}$. In their paper, Fitzpatrick et al. establish the first results on the deterministic version of the game of Zombies and Survivor [20]. They provide an example showing that if a graph is cop-win, then it is not necessarily zombie-win. They provide a sufficient condition for a graph to be zombie-win. They also establish several results about the zombie number of the Cartesian product of graphs.

The main aspect that makes different variants of these pursuit-evasion problems quite challenging is the fact that the cop number and zombie number is not a monotonic property with respect to subgraphs. For example, both the cop number and zombie number of a clique is 1 but the cop number and zombie number of a cycle on more than 3 vertices is 2.

### 1.1 Contributions

In this paper, we first consider the deterministic version of Zombies and Survivor. We then turn our attention to a deterministic variant which we call *Lazy Zombies and Survivor*. In this variant, each zombie is allowed to stay still whenever they want. The *lazy zombie number* $z_L(G)$ of a graph $G$ is therefore defined as the minimum number of lazy zombies required to catch the survivor on $G$. The *universal lazy zombie number* of a graph $G$, denoted $u_L(G)$, denotes the minimum number of lazy zombies required to catch the survivor on $G$, when the starting positions of the lazy zombies are chosen by an adversary.

We show that there exist outerplanar graphs and 2-connected outerplanar graphs $G$ with $n$ vertices such that $\frac{z(G)}{c(G)} = \Omega(n)$. This improves upon a result of Bartier et al. who showed that this ratio is $\frac{z(G)}{c(G)} = \Omega(\log n)$ for outerplanar graphs [7]. We also show that there exist maximum-degree-3 outerplanar graphs $G$ such that $\frac{z(G)}{c(G)} = \Omega\left(\frac{n}{\log(n)}\right)$ and there exist simple polygons whose visibility graph $G$ is such that $\frac{z(G)}{c(G)} = \Omega(n)$.

Then, we establish that lazy zombies are more powerful than plain zombies and less powerful than cops. Indeed, we prove that 2 lazy zombies always win (in less than $2n$ rounds) on outerplanar graphs. However, we show that there exists graphs of treewidth 2 that require 3 lazy zombies, whereas 2 cops are sufficient. We then show that $k$ lazy zombies win after $O(n^{2k})$ rounds on graphs with treedepth $k$. Finally, we highlight a few implications stemming from this upper bound such as $(k+1)\log n$ lazy zombies win on graphs with treewidth $k$, $O(\sqrt{n})$ lazy zombies are always sufficient to win on planar graphs, $O(\sqrt{gn})$ lazy zombies win on graphs with genus $g$ and $O(h\sqrt{hn})$ lazy zombies win on all connected graphs $G$ with any excluded $h$-vertex minor $H$. Our upper bounds on lazy zombie numbers still hold for universal lazy zombies.

| | zombies | (universal) lazy zombies | cops |
|---:|:---:|:---:|:---:|
| outerplanar | $\Theta(n)$ (Thm. 1) | 2 (Thm. 6) | 2 ([16]) |
| planar | $\Theta(n)$ (Thm. 1) | $O(\sqrt{n})$ (Cor. 4) | 3 ([1]) |
| genus $g$ | $\Theta(n)$ (Thm. 1) | $O(\sqrt{gn})$ (Cor. 5) | $3 + \lfloor \frac{3g}{2} \rfloor$ ([40]) |
| treedepth $k$ | $\Theta(n)$ (Thm. 1) | $k$ (Cor. 2 and Thm. 8) | $(k/2) + 1$ ([27]) |
| treewidth $k$ | $\Theta(n)$ (Thm. 1) | $(k+1)\log n$ (Cor. 2 and Lem. 2) | $(k/2) + 1$ ([27]) |
| $h$-vertex excluded minor | $\Theta(n)$ (Thm. 1) | $O(h\sqrt{hn})$ (Thm. 6) | $\frac{1}{2}(h-1)(h-2)$ ([4]) |

Table 1: Summary of zombie, (universal) lazy zombie, and cop numbers. The upper bound on the cop number for $h$-vertex excluded minors is for *connected* excluded minors.

These results are summarized in Table 1.

## 2 Preliminaries and Notation

Let $G$ be a simple and undirected graph with vertex set $V(G)$ and edge set $E(G)$. Throughout this paper $n$ will be used for $|V(G)|$. Given a subset $S \subseteq V(G)$, we denote the graph induced on $S$ as $G[S]$. Given two vertices $u$ and $v$ in $V(G)$, we denote the shortest path from $u$ to $v$ as $\pi_G(u,v)$ and its length as $|\pi_G(u,v)|$[1]. If there is no path from $u$ to $v$, then the length of the path is infinite. The diameter of a connected graph $G$, denoted as $\delta(G)$, is $\max\{|\pi_G(u,v)| : u,v \in V(G)\}$.

We denote the cop number of $G$ by $c(G)$, the zombie number of $G$ by $z(G)$ and the lazy zombie number of $G$ by $z_L(G)$. If $G$ is disconnected, then the cop (zombie, lazy zombie) number of $G$ is simply equal to the sum of the cop (zombie, lazy zombie) numbers of its connected components. Therefore, in this paper, we assume that $G$ is connected. Observe that if $k$ cops can win on $G$, then they can win wherever they start. Indeed, if their strategy involves a specific starting position, then wherever they start, they can first reach this starting position before applying their strategy. This observation does not apply to zombies. Since zombies have to follow shortest paths, they might not be able to reach a strategic starting position. In their paper, Fitzpatrick et al. [20, Figure 5] provide an example of a graph with zombie number 1, where the zombie has to start on a specific vertex or otherwise it will lose. Therefore, we define the *universal zombie number* $u(G)$ of a graph $G$ as the minimum number of zombies required to catch the survivor on $G$, where the starting position of the zombies is chosen by an adversary. Similarly, the *universal lazy zombie number* $u_L(G)$ of a graph $G$ is defined as the minimum number of lazy zombies required to catch the survivor on $G$, where the starting position of the lazy zombies is chosen by an adversary. Since universal (lazy) zombies have less power than (lazy) zombies, we have $z(G) \leq u(G)$ and $z_L(G) \leq u_L(G)$, and Fitzpatrick et al. provided an example of a graph $G$ where $z(G) < u(G)$ [20].

All the pursuit-evasion games we describe in this paper are made of a sequence of rounds, each of which is made of two turns. For each round $i \geq 0$, the zombies play first (zombies' turn) and then the survivor plays (survivor's turn). In round 0, during the zombies' turn, the zombies choose their starting position (or an adversary assigns one to them), and then, during the survivor's turn, the survivor chooses its starting position. Then the zombies move (or wait if the version of the game allows them to) and the survivor moves (or wait if they decide to) in the subsequent rounds.

---
[1] The length of the shortest path is the number of edges on the path.

When stating results about (universal) (lazy) zombie number, we sometimes use asymptotic notation. In this notation, lower-bound asymptotics ($\Omega$ or $\Theta$) will be a lower bound on the maximum-valued graph of the given class. For instance, "for outerplanar graphs $G$, $z(G) \in \Theta(n)$" means not only is the zombie number at most some constant times $n$ for every outerplanar graph, but also the zombie number is at least some constant times $n$ for some outerplanar graph family.

Throughout the paper, unless specified, the base of the logarithmic function $\log(\cdot)$ is 2.

## 3 Linear bound on zombie number

In their paper, Fitzpatrick et al. ask how large the ratio $\frac{z(G)}{c(G)}$ can be [20, Question 19]; they note that they have not observed any graph with a ratio that exceeds 2. Here we show that this ratio can be infinite and of size $\Omega(n)$, and we show this even for outerplanar graphs of fixed radius. In independent work, Bartier et al. showed that this ratio can be infinite and of size $\Omega(\log n)$, for outerplanar graphs [7].

**Theorem 1.** *Let $k \geq 2$ be an integer. Then there is a connected outerplanar graph $G_k$ with $23k + 1$ vertices that requires at least $k$ zombies.*

*Proof.* Let $H$ be the 23-vertex graph shown in Figure 1a. $H$ has two distinguished vertices $s$ and $t$. To form the graph $G_k$, first take $k$ disjoint copies $H_1, H_2, \ldots H_k$ of $H$, with $H_i$ having distinguished vertices $s_i$ and $t_i$; to this add a vertex $c$ that is connected to each $s_i$ and each $t_i$ (see Figure 1b).

By construction, $G_k$ has $23k + 1$ vertices. Suppose that $k - 1$ or fewer zombies play on $G_k$. This mean that in round 0 there is some copy $H_i$ of $H$ that contains no zombie; the survivor chooses the vertex adjacent to $s_i$ in $H_i$ as its starting position, and will stay in $H_i$ forever.

Each zombie will therefore first take a shortest path to $c$ from wherever it starts, as this is the only way to get to $H_i$. Consider the time unit on which a zombie reaches $c$; this can be anywhere from 0 to 9, as the zombie could start on $c$, and 9 is the radius of $G_k$ (and $c$ is the center). See Figure 1c. At one time unit later, the zombie will move to $s_i$ or $t_i$, whichever is closer to the survivor. We call this time the *arrival time* of the zombie; arrival times are all between 1 and 10, inclusive.

The survivor's strategy is to walk in $H_i$ away from $s_i$ until it reaches a vertex of degree three. At this point they walk along the 13-cycle in $H_i$, starting in the direction of the vertex of degree two. They continue walking this cycle forever.

The zombies will arrive on $s_i$ if their arrival time is at most five, as this is the closest vertex to the survivor at this time. These zombies will follow the survivor's path. If a zombie has arrival time six or more, it arrives on $t_i$. These zombies will pursue the survivor by first walking to the 13-cycle and then following the survivor around it.

Several steps of the chase are illustrated in Figure 2. In the figures, the survivor is at the green vertex, the zombies are at the red vertices, and the number(s) next to a red vertex indicate the arrival time of the zombies on that vertex. The illustrations end when the survivor has all zombies on the five vertices behind them on the 13-cycle. The survivor keeps walking the cycle with the five

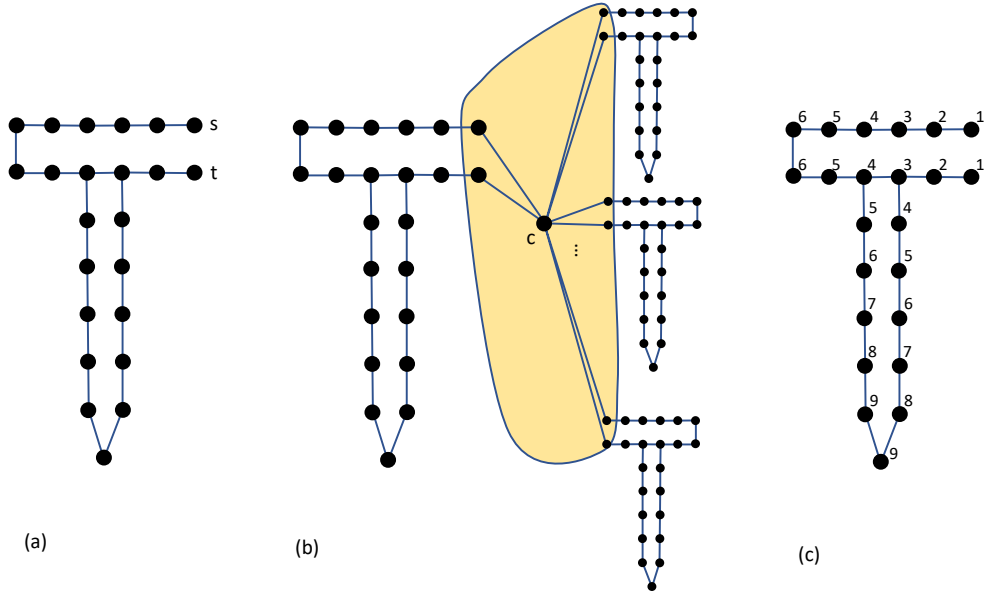

Figure 1: Construction of an outerplanar graph requiring a linear number of zombies. (a) the component graph $H$. (b) connecting components into the graph $G_k$. (c) each vertex of some $H_j$ labelled with its distance to $c$.

vertices of zombies behind them forever.

Therefore, since $k-1$ zombies are insufficient to capture the survivor on $G_k$, at least $k$ are required and the lemma is proved. □

Note that a linear number of zombies always suffices for a graph, as we could use $n$ zombies and initially place one on each vertex (or perhaps leave one free for the survivor). Thus we have shown that for general or for outerplanar graphs, $z(G) \in \Theta(n)$. Since the copnumber for outerplanar graphs is at most two 2, the ratio $\frac{z(G_k)}{c(G_k)}$ is $\frac{k}{2} = \frac{n-1}{46} \in \Theta(n)$.

Modifications of the construction in the proof of Theorem 1 will work for other graph classes.

**Theorem 2.** *Let $k \geq 2$ be an integer. There is a 2-connected outerplanar graph $G_k$ with $30k+1$ vertices that requires at least $k$ zombies.*

*Proof.* We start with a different component graph $H$ shown in Figure 3a. As in the proof of Theorem 1, we connect each $s_i$ and $t_i$ to a new center vertex $c$. To make this 2-connected, we connect $t_i$ to $s_{i+1}$ for $1 \leq i \leq k-1$.

The proof is similar to the previous one. The main difference is that at time 7, zombies can arrive at both $s_i$ (from $c$ or $t_{i-1}$) and $t_i$ (from $c$ or $s_{i+1}$). □

**Theorem 3.** *Let $k \geq 2$ be an integer. There is a connected graph $G_k$ with $15k$ vertices that requires at least $k$ zombies.*

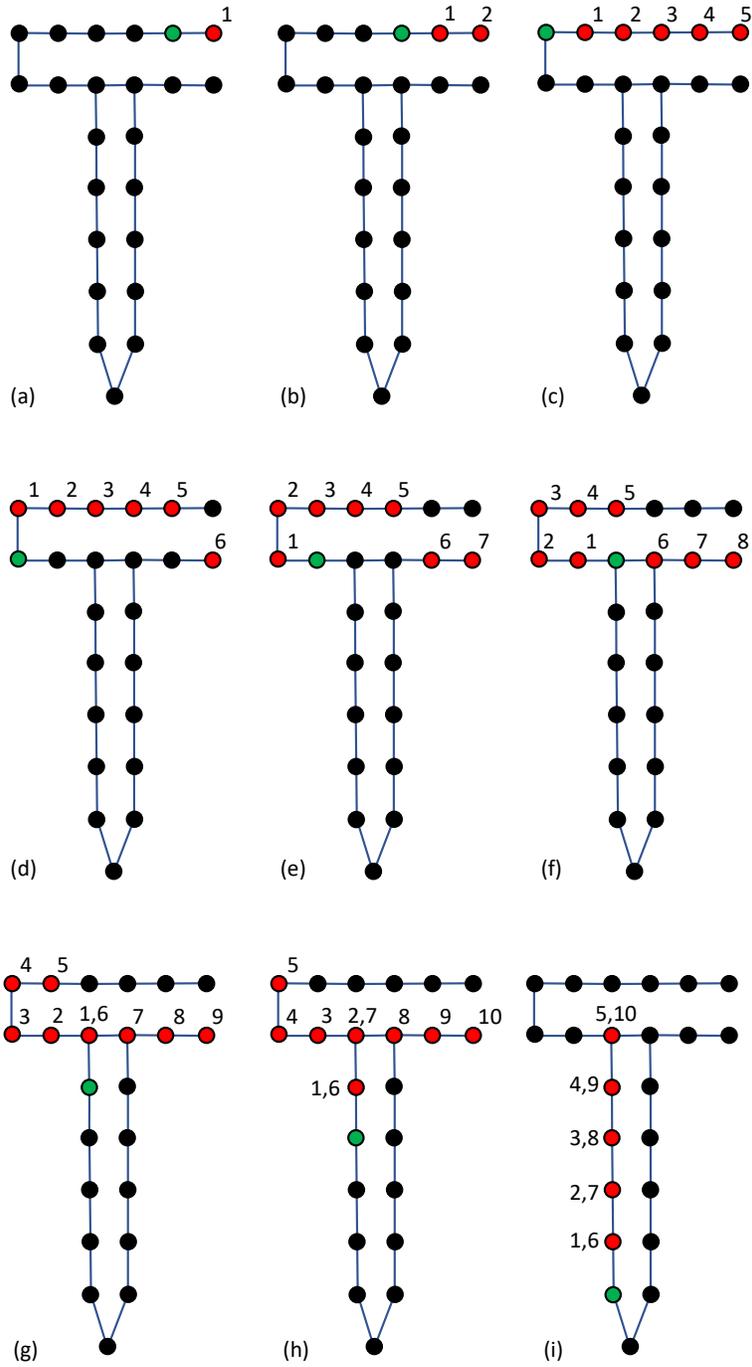

Figure 2: Steps in the survivor strategy for $k-1$ zombies in $G_k$. (a) turn 1, (b) turn 2, (c) turn 5, (d) turn 6, (e) turn 7, (f) turn 8, (g) turn 9, (h) turn 10, (i) turn 13.

*Proof.* Use $k$ copies of the 15-vertex component $H$ shown in 3b, and connect all $s_i$ and $t_i$ into a clique. The rest of the proof is similar to that of Theorem 1, with the arrival times between one and

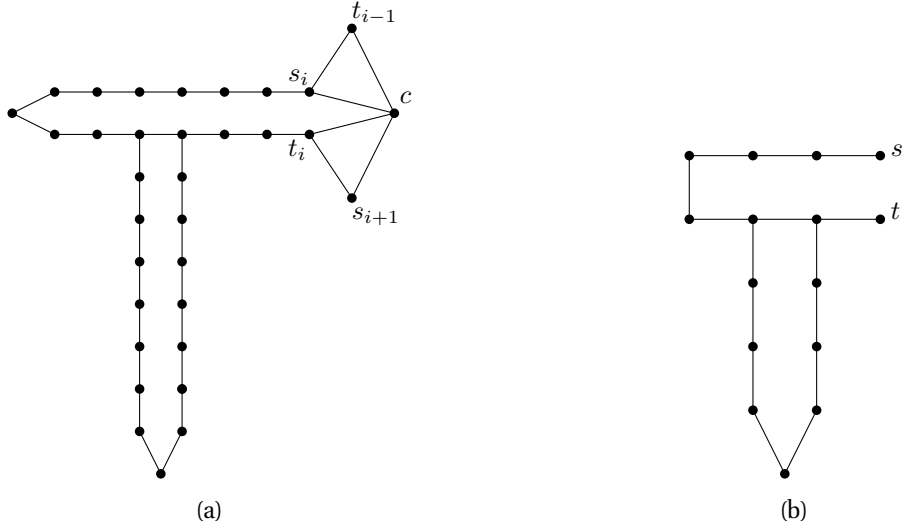

Figure 3: (a) component graph $H$ used in the proof of Theorem 2. (b) component graph $H$ used in the proof of Theorem 3.

six. □

**Theorem 4.** *Let $k \geq 2$ be an integer. There is a maximum-degree-3 connected outerplanar graph $G_k$ with at most $24k + 16k\lceil \log k \rceil - 1$ vertices that requires at least $k$ zombies.*

*Proof.* Consider replacing the center vertex $c$ in the construction of the proof of Theorem 1 with a rooted, minimum-height binary (degree 3, except the root) tree $T$ having $k$ leaves. Each leaf $l_i$ of the tree will connect to a single $s_i$ and associated $t_i$. Since there are $k$ different copies of the component graph $H$, the tree will need $k$ leaves, which can be done with a tree of height $\lceil \log k \rceil$.

This extends the range of arrival times for zombies. Suppose the survivor has chosen $H_i$ for its initial location. A zombie could have started on $l_i$; this zombie would have an arrival time of 1. The latest arrival time happens when a zombie is as deep as possible in some $H_j$ (with $j \neq i$), and $l_i$ and $l_j$ are at a distance of twice the height of the tree $T$. So if $a$ is the maximum possible arrival time of a zombie,

$$a = \max_{v \in H_j} d(v, l_j) + 2\lceil \log k \rceil + 1, \qquad (*)$$

where $d$ is the graph-theoretical distance.

Since the $H$ chosen in the proof of Theorem 1 is made for arrival times up to 10, we have to modify it for this situation. Refer to Figure 4. We will design our new $H$ to handle arrival times of up to $a^*$, where $a^*$ is $4\lceil \frac{a-2}{4} \rceil + 2$, the smallest number equivalent to 2 mod 4 that is at least as large as $a$. $H$ will have a shortest $s-t$ path of length $a^* + 1$, with a cycle of length $a^* + 3$ sharing the $\lceil \frac{a^*}{4} \rceil$ edge from $t$ with it. This $H$ has $2a^* + 3$ vertices. The maximum distance to escape a copy $H_i$ of $H$ to $l_i$ is $3\lceil \frac{a^*}{4} \rceil = 3\frac{a^*+2}{4}$. Substituting this in for the max in (*), we get

$$a = \frac{3a^* + 6}{4} + 2\lceil \log k \rceil + 1.$$

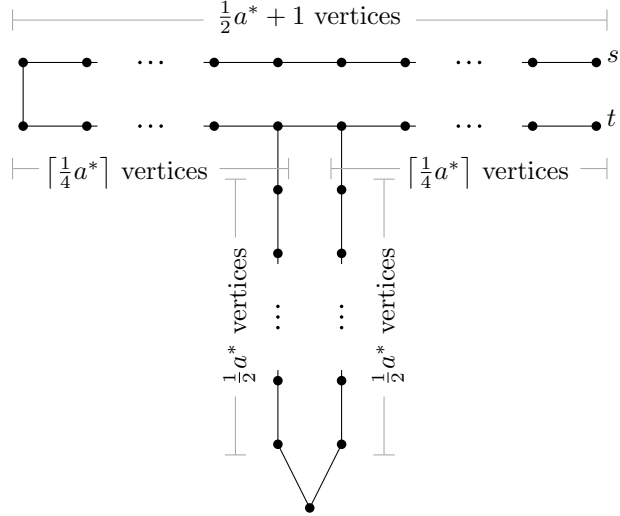

Figure 4: Component graph $H$ used in the proof of Theorem 4.

Since $a \leq a^*$, we must satisfy
$$\frac{3a^* + 6}{4} + 2\lceil \log k \rceil + 1 \leq a^*.$$

Simplifying, we arrive at
$$a^* \geq 8\lceil \log k \rceil + 10.$$

Furthermore, we know that $a^* \equiv 2 \pmod{4}$, and should be as small as possible. Thus
$$a^* = 8\lceil \log k \rceil + 10.$$

The number of vertices of the entire construction is $k-1$ for the tree T, plus $k$ times $2a^* + 3$ for the $k$ copies of $H$. This is at most $k - 1 + 2k(8\lceil \log k \rceil + 10) + 3k = 24k + 16k\lceil \log k \rceil - 1$ vertices.

Proof that the survivor cannot be captured with fewer than $k$ zombies is analogous to the proof in Theorem 1, with the survivor starting adjacent to $s$, walking the $s - t$ path until it reaches a degree-3 vertex, and then walking the cycle. The graph is maximum degree 3 by construction. □

The previous theorem gives us $n \in O(k \log k)$, or $n \leq ck \log k$ for some constant $c$. Hence, we have $\frac{n}{\log n} \leq \frac{n}{\log k} \leq ck$, or $k \geq \frac{n}{c \log n}$. Since $k$ zombies are required, this gives us a lower bound of $\Omega(\frac{n}{\log n})$ on the zombie number of bounded-degree graphs.

### 3.1 Polygon visibility graphs

A *polygonal chain* is a finite sequence $V$ of points $v_1, v_2, \ldots, v_n$ in $\mathbf{R}^2$ (called *vertices*) along with the line segments $v_1 v_2, v_2 v_3, \ldots v_{n-1} v_n$ (called *edges*). A polygonal chain is called *closed* if $v_1 = v_n$, and *simple* if no two edges intersect except consecutive edges intersecting at their common vertex. A closed simple polygonal chain divides the plane into a finite *interior* and infinite *exterior*. A *simple polygon*, or simply *polygon*, is a closed simple polygonal chain along with its interior.

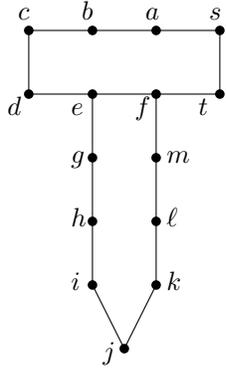
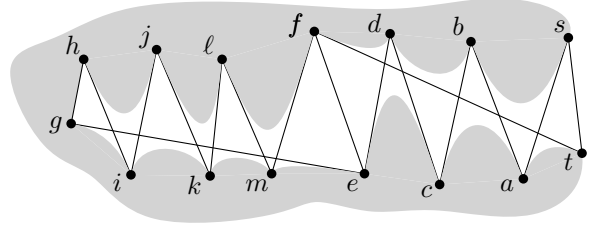

(a) A graph fragment for proving lower bounds.

(b) An embedding of the graph fragment in the visibility graph of a polygon fragment.

Figure 5: Creation of a polygon fragment to prove a lower bound.

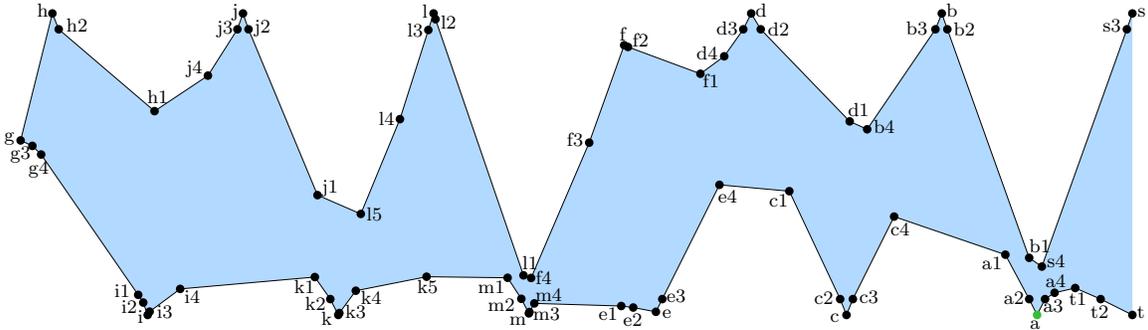

Figure 6: The polygon fragment Q. Note that at this scale, the dots for the vertex-pair $(i, i3)$ merge into what can appear to be a single dot. This also happens for $(k, k3)$ and $(m, m3)$.

The *visibility graph* of a simple polygon $P$ [17, 32] is a graph $G$ where $V(G)$ is the set of vertices of the polygon, and $E(G) = \{(v_i, v_j) |$ the segment $v_i v_j$ does not intersect the exterior of $P\}$. In particular, this means that every edge of the polygon is an edge of the visibility graph, but the visibility graph has other edges corresponding to segments that traverse the interior and possibly boundary of $P$. Visibility graphs are of interest in discrete geometry and have applications, for instance, in motion planning and shape analysis [15, 17].

Here we show that there is a linear bound on the zombie number of visibility graphs. The proof is messy in that it involves a large polygonal chain with relatively precise vertex locations in order to get a visibility graph with the desired properties. However, it is inspired by the proof of Theorem 1. Consider the graph fragment in Figure 5a. One way to embed this graph inside a polygon visibility graph is sketched in Figure 5b. Complicating matters is that we can't get the required non-edges without placing vertices inbetween those shown in Figure 5b. Once those non-visibilities are worked out, we get the polygon fragment $Q$ shown in Figure 6. For reproducibility, the exact vertex locations of $Q$ are given in Appendix A (refer to Table 2).

To form a polygon whose visibility graph requires at least $k$ zombies, we will connect $k$ copies of $Q$, denoted $Q_1, Q_2, \ldots Q_k$, placed in a geometric configuration where the only vertices of $Q_i$ visible to $Q_j$ are $s_i$ (the copy of $s$ in $Q_i$) and $t_i$ (similar). Thus $S$, the collection of all $s_i$'s and $t_i$'s, will form a clique in the visibility graph. This is done by taking a small sliver of a circular arc and placing the $2k$ vertices of $S$ evenly along it, as shown in Figure 7. If the sliver is small enough, any vertices inside $Q_i \cap \overline{S}$ will see only the relative interior of the polygon edge $t_k s_1$. (Another method of ensuring this is to scale each $Q_i$ up in its $x$-coordinate, effectively pushing interior vertices away from $s_i t_i$. Such a scale does not affect the visibility graph of $Q_i$.)

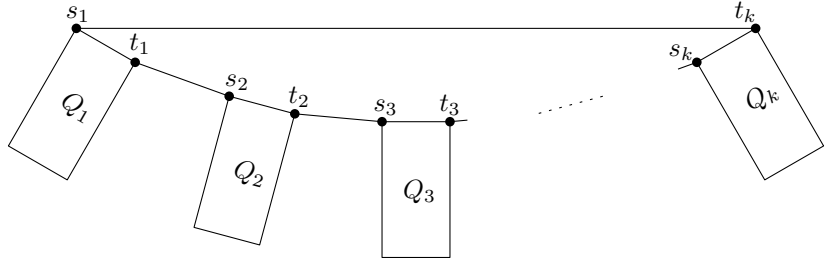

Figure 7: Connecting the $Q_i$'s into a polygon.

The proof that this polygon has zombie number $k$ now roughly follows that of Theorem 1. Suppose that less than $k$ suffices. Then, in the zombies' initial placement, there will be one copy of $Q$, say $Q_i$, that has no zombies in it. Start with the survivor on vertex $a$ of $Q_i$ (refer to Figure 5a for the lebelling of the vertices).

$Q$ is constructed so that it takes at most 6 turns for a zombie to leave the fragment (or 5 turns at most to get to $s$ or $t$). This means that if the survivor stays in $Q_i$ for 6 turns (it will), all zombies will have arrived in $Q_i$. The survivor's strategy will be to walk from $a$ to $b$, to $c$, to $d$, $e$, $g$, $h$, $i$, $j$, $k$, $l$, $m$, and $f$. Unlike in Theorem 1, the survivor cannot now loop back to $e$ (they might be caught by a zombie) but must instead move to $t = t_i$. Once at $t_i$, the survivor chooses some $s_j$ where $i \neq j$, and moves there. Next they can move to $a_j$ and start the same walk in $Q_j$ as it did in $Q_i$. It may continue in this way *ad infinitum*. In Appendix A (refer to Figures 15 to 22), we show a few of the different steps in the survivor's initial walk through $Q_i$.

Since $Q$ has 69 vertices, we have shown the following.

**Theorem 5.** *Let $k \geq 2$ be an integer. Then there is a polygon $P_k$ with $69k$ vertices whose visibility graph requires at least $k$ zombies.*

A linear number of zombies will always work (e.g. $n/3$ of them starting on every third vertex), so the maximum zombie number of the visibility graph of a polygon with $n$ vertices is $\Theta(n)$.

There is a related problem that asks for the zombie number of a *point-visibility graph* of a polygon, which is the infinite graph $G$ where $V(G)$ is taken to be the *points* of the polygon, not simply the vertices. Edges are then defined as in the visibility graph. This problem involves more geometry than the other problems we have studied. Here it is not clear that there are polygons with a point-visibility graph zombie number higher than one.

## 4 The Lazy Zombie Number of Outerplanar Graphs is 2

In the previous section, we showed that $\Omega(n)$ zombies are sometimes necessary to catch a survivor on an outerplanar graph. In this section, we show that 2 lazy zombies are always sufficient to catch the survivor on outerplanar graphs. Observe that two lazy zombies are sometimes necessary to catch a survivor on an outerplanar graph since a single lazy zombie cannot win on a 4-cycle.

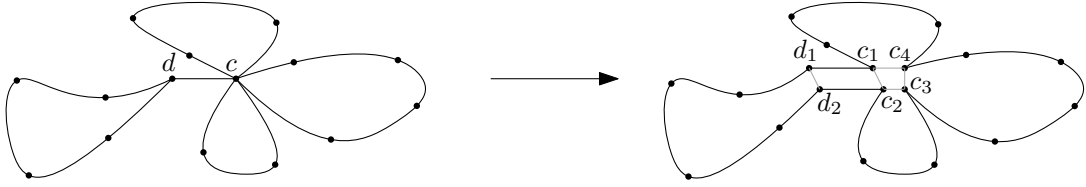

Figure 8: We duplicate each cut edge and we add a special null chord between every pair of consecutive appearances of each cut vertex on the circuit.

**Theorem 6.** *Let G be a connected outerplanar graph. Then $z_L(G) = 2$ and 2 lazy zombies can catch the survivor in less than $2n$ rounds.*

The proof of Theorem 6 has some points in common with the proof that the cop number of outerplanar graphs is 2 [16]. However, in our case, every time a lazy zombie is moving, we need to make sure that it does so along a shortest path to the survivor.

*Proof.* If $G$ is a tree or a cycle, then 2 lazy zombies can win on $G$ since 2 plain zombies can win on $G$ [20]. Otherwise, the boundary of the outerface is a circuit where cut vertices can appear several times and cut edges appear twice [18]. Any edge not on the circuit is called a *chord*. We duplicate each cut edge and we add a special *null chord* between every pair of consecutive appearances of each cut vertex on the circuit (refer to Figure 8).

In the following, we will have one lazy zombie (denoted $z_1$) be *stationary*, initially on one end of a chord or null chord $b_i b_j$. This lazy zombie will capture the survivor if the survivor moves to $b_i$ or $b_j$ but otherwise will not move. We refer to the vertices of $G$ where the survivor is known to be restricted as the *survivor territory*. The remaining vertices of $G$ are the *zombie territory*. As the game evolves, the survivor territory (and hence the zombie territory) will change. At the beginning, the circuit is divided into two or more connected components by the removal of the vertices $b_i$ and $b_j$; the survivor is in one of these connected components. The vertices of this connected component is called the initial survivor territory. The remaining vertices of $G$ are the initial zombie territory. If given a choice of chords for $b_i b_j$, we choose the one that limits the survivor territory to the fewest vertices possible. Therefore, by outerplanarity, there is no chord or null chord from the survivor territory to the zombie territory.

The other lazy zombie (denoted $z_2$), the *advancing* lazy zombie, will start either on the chord $b_i b_j$ or in the zombie territory. On each turn, it advances toward the survivor. If it starts in the zombie territory, then at some point it will be either at $b_i$ or $b_j$. Without loss of generality, we can assume that $z_2$ is at $b_j$ and there is a shortest path from $b_j$ to the survivor whose first edge has an endpoint in the survivor territory. Otherwise, $z_2$ is at $b_j$. If all shortest paths to the survivor go through $b_i$, then $z_2$ moves to $b_i$. Then, in the next turn, it cannot be the case that all shortest paths to the survivor go through $b_j$. Hence, $z_2$ is at $b_i$ and there is a shortest path from $b_i$ to the survivor whose first edge has an endpoint in the survivor territory. We are in a symmetric situation ($z_1$ is on one end of $b_i b_j$ and $z_2$ is at $b_i$ with a shortest path to the survivor entering the survivor territory).

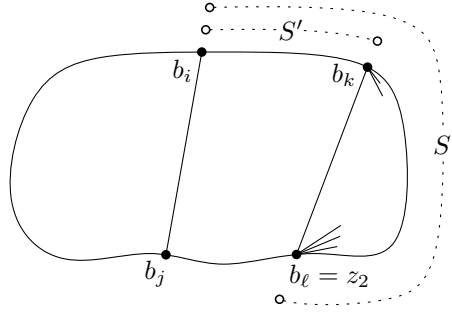

Figure 9: The advancing lazy zombie is on a vertex $b_\ell$ with a chord (or null chord) whose other end is in $S$.

We use the term *interval* to refer to the vertices of a walk that is a connected component of the circuit which is the boundary of the outerface. We will maintain the invariant that (refer to Figure 9)

1. the survivor territory corresponds to the clockwise interval $S$ between $b_i$ (exclusive) and $z_2$ (exclusive),

2. there is no chord or null chord from $S$ to the clockwise interval between $z_2$ (exclusive) and $b_j$

3. and there is a shortest path from $z_2$ to the survivor whose first edge has an endpoint in $S$.

The advancing lazy zombie moves as follows:

- If $z_2$ is adjacent to the survivor, it captures the survivor.

- If $z_2$ is on a vertex with no chord or null chord whose other end is in $S$, it moves to its counterclockwise neighbour on the circuit (this is possible by our invariant). Afterwards, our invariant is maintained and the survivor territory decreased by one vertex.

- If $z_2$ is on a vertex $b_\ell$ with a chord (or null chord) whose other end is in $S$, let $b_k b_\ell$ be the chord closest to and clockwise of $b_i$ (refer to Figure 9).

    – If the survivor is in the clockwise interval between $b_k$ and $b_\ell$, then $z_2$ becomes the stationary lazy zombie holding $b_k b_\ell$ shut, and $z_1$ becomes the advancing lazy zombie. Then $z_1$ moves until it reaches $b_k b_\ell$. When it does, after at most one extra move by $z_1$, our invariant is restored and the survivor territory decreased by at least one vertex.

    – If the survivor is in the clockwise interval between $b_i$ and $b_k$, then $z_2$ advances to $b_k$. Let $S'$ be the clockwise interval between $b_i$ (exclusive) and $b_k$ (exclusive). By our invariant, there is no chord or null chord from $S'$ to the clockwise interval between $b_\ell$ (exclusive) and $b_j$. By our choice of chord, there is no chord from $S'$ to $b_\ell$. Moreover, by outerplanarity, there is no chord or null chord from $S'$ to the clockwise interval between $b_k$ (exclusive) and $b_\ell$ (exclusive). Thus our invariant is maintained with $S' \subsetneq S$ as the survivor territory.

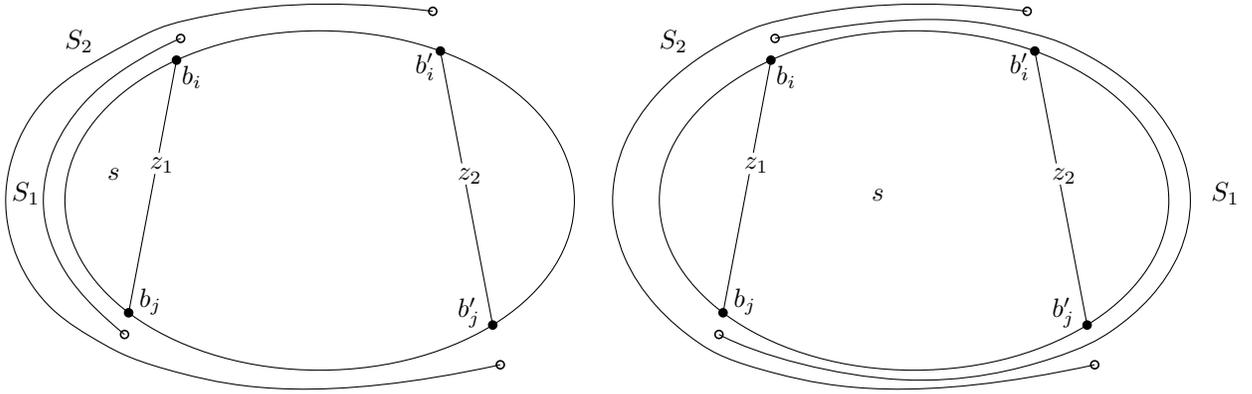

Figure 10: The survivor territories $S_1$ and $S_2$ defined by $z_1$ and $z_2$, depending on the position of $s$. If $s$ is in the clockwise interval between $b'_i$ (exclusive) and $b'_j$ (exclusive), we get a configuration that is symmetric to the one on the left.

The advancing lazy zombie's walk ends with either a capture or a reduction (in the number of vertices) in the survivor territory. Eventually, the survivor's territory is reduced to a single path and then the survivor will be captured.

Assume the lazy zombies follow the strategy described in this proof. Then, in each round, exactly one lazy zombie is moving. Moreover, for all vertices $v \in V$, each lazy zombie will visit $v$ at most once. Therefore, the capture happens in less than $2n$ rounds. □

**Corollary 1.** *Let $G$ be a connected outerplanar graph. Then $u_L(G) = 2$.*

*Proof.* Let some adversary choose the starting positions $z_1$ and $z_2$ for two lazy zombies. If $G$ is a tree, then a single lazy zombie will win by following the unique path towards the survivor. If $G$ is a cycle, then one lazy zombie stays stationary and the other one chases the survivor (always in the same direction). The survivor will eventually get squeezed between the two lazy zombies. Otherwise, as in the proof of Theorem 6, the boundary of the outerface is a circuit where cut vertices can appear several times and cut edges appear twice [18]. Any edge not on the circuit is called a *chord*. We duplicate each cut edge and we add a special *null chord* between every pair of consecutive appearances of each cut vertex on the circuit (refer to Figure 8).

Each lazy zombie $z_\ell$ ($\ell = 1, 2$) will first proceed as follows: $z_\ell$ chases the survivor $s$ until it reaches a chord. As soon as $z_\ell$ reaches a chord, it stops and waits for the other lazy zombie to reach a chord as well. If $z_1$ and $z_2$ both reached the same chord, then we apply the strategy described in the proof of Theorem 6. Otherwise, let $b_i b_j$ be the chord protected by $z_1$ and $b'_i b'_j$ be the chord protected by $z_2$. Without loss of generality, assume that $b_i$, $b'_i$, $b'_j$ and $b_j$ appear in this order, clockwise, around the outerface (refer to Figure 10). The chord that is covered by $z_\ell$ splits the graph into a zombie territory $Z_\ell$ and a survivor territory $S_\ell$ (refer to the proof of Theorem 6). If given a choice of chords, $z_\ell$ chooses the one that limits $S_\ell$ to the fewest vertices possible. Therefore, we have $s \in S_1 \cap S_2$. If $z_1 \in Z_2$ or $z_2 \in Z_1$ (refer to Figure 10 (left)), then we apply the strategy described in the proof of Theorem 6.

Otherwise, we have $z_1 \in S_2$ and $z_2 \in S_1$ (refer to Figure 10 (right)). From there, we make

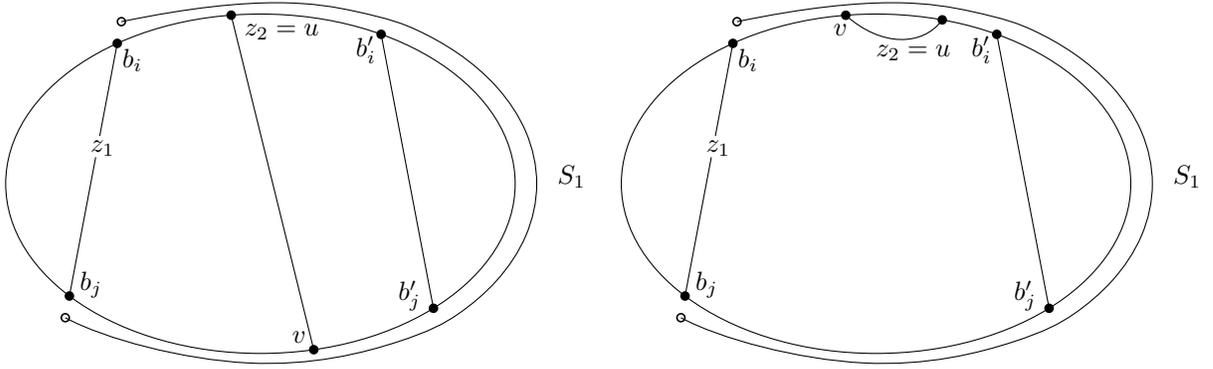

Figure 11: Two possible locations for the edge $uv$. Case (1) is depicted on the left. If $s$ is in the clockwise interval between $u$ (exclusive) and $v$ (exclusive), then $z_1 \in Z_2'$. Case (2) is depicted on the right. If $s$ is in the clockwise interval between $v$ (exclusive) and $u$ (exclusive), then $z_1 \in Z_2'$.

$z_1$ the stationary lazy zombie and $z_2$ starts chasing $s$. Without loss of generality, assume that $z_2$ starts the chase by moving counterclockwise. Observe that, during this chase, if at a given round $z_2$ travels along an edge $uv$ from $u$ to $v$, then in the next round, it cannot be the case that all shortest paths from $z_2$ to the survivor go from $v$ to $u$. So as the chase is going on, $z_2$ will keep moving counterclockwise. Here is what happens when $z_2$ encounters a chord $uv$. Assume $z_2$ is on $u$. The chord $uv$ splits the graph into a new zombie territory $Z_2'$ and a new survivor territory $S_2'$. We consider two cases: (1) $u$ is in the clockwise interval between $b_i$ and $b_i'$ (exclusive), and $v$ is in the clockwise interval between $b_j'$ and $b_j$, or (2) both $u$ and $v$ are in the clockwise interval between $b_i$ and $b_i'$ (exclusive). (refer to Figure 11)

(1) We consider two subcases: (a) $z_1 \in Z_2'$ or (b) $z_1 \notin Z_2'$.

   (a) If $z_1 \in Z_2'$, then we apply the strategy described in the proof of Theorem 6.
   
   (b) Otherwise, we are in the same situation again, that is, $s \in S_1 \cap S_2'$, $z_1 \in S_2'$ and $z_2 \in S_1$. However, the size of $S_1 \cap S_2'$ is now smaller than $S_1 \cap S_2$ by at least one vertex. So $z_2$ restarts the chase for $s$.

(2) We consider two subcases: (a) $z_1 \in Z_2'$ or (b) $z_1 \notin Z_2'$.

   (a) If $z_1 \in Z_2'$, then we apply the strategy described in the proof of Theorem 6.
   
   (b) Otherwise, $z_2$ moves to $v$. From there, $z_2$ keeps chasing the survivor by moving counterclockwise. If it falls into Case (2)(b) again, observe that $z_2$ is now closer to $b_i$. So falling back into Case (2)(b) can only happen a finite number of times in a row. Eventually, $z_2$ will fall into Case (1) or Case (2)(a).

When $z_2$ gets into Case (1)(a) or Case(2)(a), then we apply the strategy described in the proof of Theorem 6 and the survivor gets caught. Every time $z_2$ gets into Case (1)(b), then we are closer to a configuration where we can apply the strategy described in the proof of Theorem 6. Every time $z_2$ gets into Case (2)(b), then after a finite number of rounds, it falls into Case (1) or Case (2)(b). In all cases, the survivor will eventually get caught. □

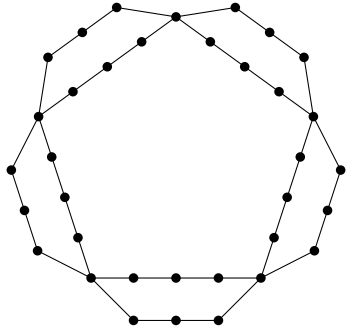

Figure 12: An example of a treewidth-2 graph $G$ with $z_L(G) > 2$.

Outerplanar graphs are a subset of the treewidth-2 graphs, but Theorem 6 cannot be generalized to treewidth-2. Figure 12 shows a graph with treewidth two that requires three lazy zombies. For two lazy zombies, the survivor strategy is based on travelling between vertices of degree four. The survivor starts at any degree-four vertex that is not adjacent to a lazy zombie, and waits at this vertex until it is adjacent to a lazy zombie. If it is adjacent to two lazy zombies, it picks a free neighbor and moves to it, continuing to travel until it reaches the next vertex of degree four. If it is adjacent to only one lazy zombie, then it chooses one of the three free neighbors that leads away from the other lazy zombie, moving to the vertex of degree four reachable through that free neighbor. Once the survivor arrives on a vertex of degree four, it waits (possibly zero turns) until it is adjacent to a zombie, and repeats the procedure.

The above example shows a distinction between the lazy zombie number and the cop number of a graph since 2 cops are sufficient for a graph of treewidth 2[27]. We will study general treewidth-$k$ graphs in Section 5.

## 5 Cut-decomposable Graphs and Lazy Zombies

In this section, we explore the relationship between lazy zombie numbers and various graph parameters. We first define some of these graph parameters and some notation most of which appears in Diestel [18]. We then present the general approach, and finally, we outline some of the consequences of our approach.

Let $T$ be a tree rooted at a vertex $r$. For a vertex $v \in V(T)$, we denote the unique path from $v$ to $r$ as $\pi_T(v)$. The depth of $v$, $d_T(v) := |\pi_T(v)|$, is the length of the path from $v$ to $r$. If a vertex $u \in V(T)$ is in $\pi_T(v)$ then $u$ is an ancestor of $v$ and $v$ is a descendant of $u$. Note that $v$ is an ancestor and descendant of itself. The height of $T$ is defined as $H_T := \max\{|\pi_T(v)| : v \in V(T)\}$. The subtree of $T$ rooted at $v$ is denoted as $\Lambda_T(v)$. The height of $\Lambda_T(v)$ is defined as $H_T(v) := \max\{|\pi_T(x)| - |\pi_T(v)| : x \in \Lambda_T(v)\}$. The closure of $T$, denoted as $clos(T)$, is $T \cup \{uv : u$ is an ancestor of $v$ in $T\}$. The *treedepth* of a connected graph $G$, which we denote as $td(G)$, is 1 plus the minimum height $H_T$ over all trees $T$ defined on $V(G)$ such that $G \subseteq clos(T)$.

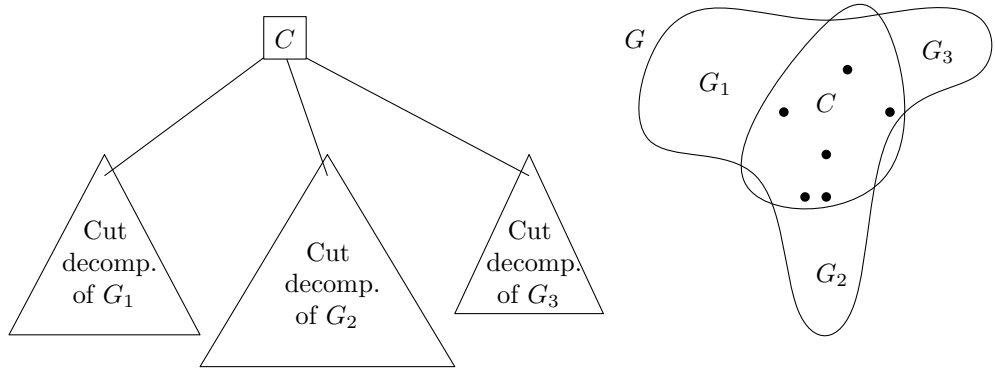

Figure 13: A cut decomposition of $G$.

The treedepth of a graph has a simple recursive definition [34].

$$td(G) = \begin{cases} 1 & \text{if } |V(G)| = 1 \\ 1 + \min_{v \in V(G)} td(G - v) & \text{if } G \text{ is connected} \\ \max_{i=1}^{p} td(G_i) & \text{otherwise, where each } G_i \text{ is a connected component} \end{cases}$$

A *tree decomposition* of a graph $G$ is a pair $(T, \mathcal{B})$ where $T$ is a tree and $\mathcal{B} = \{B_x \subseteq V(G) : x \in V(T)\}$, where each $B_x$ is a subset of $V(G)$ indexed by the nodes of $T$. The set $B_x$ is sometimes referred to as a *bag*. The following properties must be satisfied:

- For every $v \in V(G)$, the set $\{x \in V(T) : v \in B_x\}$ induces a non-empty subtree of $T$.
- For every $uv \in E(G)$, $\exists x \in V(T)$ such that $u$ and $v$ are both in $B_x$.

The *width* of a tree decomposition is 1 less than the cardinality of the largest bag. The *treewidth* of a graph $G$, which we denote as $tw(G)$, is the minimum width over all tree decompositions of $G$.

A *cut decomposition* of a graph $G$ is a pair $(X, \mathcal{C})$ where $X$ is a rooted tree and $\mathcal{C} = \{C_x \subseteq V(G) : x \in V(X)\}$, where each $C_x$ is a subset of $V(G)$ indexed by the nodes of $X$. We will refer to the set $C_x$ as the *container* of $x$ to avoid confusion with a bag of a tree decomposition. We will refer to the size of the largest container as the width of the cut decomposition tree, denoted as $cdw(X)$. The following properties must be satisfied:

- For every $v \in V(G)$, there is a unique $x \in V(X)$ such that $v \in C_x$,
- For every $uv \in E(G)$, $\exists x, y \in V(X)$ such that $u \in C_x$, $v \in C_y$ and $x$ is an ancestor of $y$ in $X$.
- For each non-leaf node $y \in X$, $C_y$ is a cut-set of $G[Y]$ where $Y = \bigcup_{x \in \Lambda_X(y)} C_x$,

Throughout this section, we will assume that $(X, \mathcal{C})$ is a cut decomposition of a graph $G$ where $|V(G)| = n$. We will refer to $X$ as a *cut decomposition tree*. We will refer to the vertices of $G$ as vertices and the vertices of $X$ as nodes, in an attempt to make the distinction clear. We will use $u, v$ to refer to vertices of $G$ and $x, y$ to refer to nodes in $X$ (or nodes in a tree decomposition). For a node $y \in X$, we define the *component* of $y$ to be $G[Y]$ where $Y = \bigcup_{x \in \Lambda_X(y)} C_x$. We slightly abuse notation and refer to the component of $y$ as $G[\Lambda_X(y)]$. Intuitively, a cut decomposition tree is a decomposition of a graph by cuts where an internal node $x$ of the tree $X$ represents a cut set of the graph $G[\Lambda_X(x)]$. The container of the root of the tree contains either the entire vertex set of $G$, or the vertices of a cut set of $G$. If it contains a cut, then the children of the root recursively correspond to the different connected components of the graph that result when the cut is removed, as in Figure 13. Cut decompositions are related to both tree decompositions and treedepth. It is related to tree decompositions except that the bags in tree decompositions have different properties than the containers of cut decompositions. In treedepth, we have a tree $T$ such that $G \subseteq clos(T)$, so for an edge $uv \in E(G)$, $uv \in clos(T)$. Cut decompositions are related, since in them, for any edge $uv \in E(G)$ where $u \in C_x$ and $v \in C_y$, $xy \in clos(X)$. We will elaborate on these relationships in the sequel.

We define the *load* of a node $x$ in $X$ as:

$$load(x) = \begin{cases} |C_x| & \text{if } v \text{ is a leaf} \\ |C_x| + \max_y load(y) & \text{otherwise, where } y \text{ is a child of } x \end{cases}$$

The load of a cut decomposition is defined as the load of the root of the cut decomposition tree. We define the load of a graph $G$, denoted as $load(G)$, to be the minimum load among all cut decompositions of $G$. We will show that $load(G)$ is a sufficient number of lazy zombies to catch a survivor in $G$.

$$time(x) = \begin{cases} |C_x|(\delta(G)-1)+1 & \text{if } x \text{ is a leaf} \\ \max_y time(y)(|C_x|(\delta(G)-1)+1) & \text{otherwise, where } y \text{ is a child of } x \end{cases}$$

The time of the root is an upper bound on the number of rounds it takes the lazy zombies to capture the survivor.

In the following, each lazy zombie $z_i$ may be *assigned* to a vertex $v$ of $G$. The strategy of the lazy zombie will be the following. If the zombie is not assigned to any vertex, then on its turn to move, it remains at its current location. A zombie assigned to a vertex $v$ has the following behavior: on its turn, it moves off its current vertex $u$ to an adjacent vertex $w$ only if there exists a $w$ that is closer to both $v$ and the survivor. This is precisely where we use the power of a lazy zombie to stand still where regular zombies cannot. Because a shortest path from $z_i$'s location to $v$ has at most $\delta(G)$ edges, the survivor can encounter vertex $v$ at most $\delta(G) - 1$ times (at or after the time $z_i$ was assigned to $v$) without being immediately caught. Lazy zombies can and will be reassigned to different vertices during the game.

We proceed by induction. The following lemma establishes the basis.

**Lemma 1.** *Let $G$ be a connected graph with cut decomposition $(X, \mathcal{C})$. Suppose that the survivor is restricted to the vertices of $G[C_x]$ for some leaf $x$ in $X$. Then, $load(x)$ lazy zombies, starting from anywhere in $G$, can capture the survivor in at most $time(x)$ rounds.*

*Proof.* Since $load(x)$ for a leaf $x$ is $|C_x|$, we assign zombies $z_1, z_2, \ldots z_{|C_x|}$ each to a different vertex of $C_x$. The survivor may visit any vertex $v \in C_x$ at most $\delta(G) - 1$ times before the zombie assigned to $v$ is adjacent to it. This is because each time the survivor lands on a vertex $v$, the zombie assigned to $v$ can move closer to the survivor. Since the survivor is restricted to vertices in $C_x$, after $|C_x|(\delta(G)-1)$ rounds, every zombie will either be on or adjacent to its assigned vertex in $C_x$. Thus, the survivor will be caught in the next round giving an upper bound on the number of rounds sufficient to catch the survivor of $time(x) = |C_x|(\delta(G)-1) + 1$. □

**Theorem 7.** *Let $G$ be a connected graph with cut decomposition $(X, \mathcal{C})$. Suppose that $x$ is a node of $X$ and the survivor is restricted to the vertices in $G[\Lambda_X(x)]$. Then, $load(x)$ lazy zombies, starting from anywhere in $G$, can capture the survivor in at most $time(x)$ rounds.*

*Proof.* By induction on $H_X(x)$, the height of $\Lambda_X(x)$. The basis, $H_X(x) = 0$, i.e. when $x$ is a leaf, follows from Lemma 1.

We assume that $load(x)$ lazy zombies are sufficient to catch the survivor in $time(x)$ rounds when the survivor is restricted to $G[\Lambda_X(x)]$, where $H_X(x) \leq k$ for $k \geq 0$. We now proceed with the case when $H_X(x) = k+1$. Let $c$ be the maximum load of a child of $x$, and $d$ be the maximum time for a child of $x$. We allocate zombies $z_1, z_2, \ldots z_c$ to the children of $x$. These zombies are initially unassigned to any specific vertex but will be assigned to specific vertices depending on the survivor's moves. Note that it is not necessarily the case that we need to use this many lazy zombies, but this number is always sufficient. We assign zombies $z_{c+1}, z_{c+2}, \ldots z_{c+|C_x|}$, each to a different vertex of $C_x$, respectively. Since $load(x) = |C_x| + c$, we have a sufficient number of zombies.

The survivor may now encounter each vertex of $C_x$ at most $(\delta(G) - 1)$ times without immediately being caught in the next round, which again follows from the upper bound of $\delta(G)$ edges on any shortest path between two vertices in $G$. Before the survivor's first encounter with a vertex of $C_x$, or between successive visits of vertices in $C_x$, or after the last visit, the survivor is restricted to the vertices of the component of the subtree rooted at exactly one child $y$ of $x$. This follows from the fact that $C_x$ is a cut set for $G[\Lambda_X(x)]$. We apply the inductive hypothesis on $G[\Lambda_X(y)]$, since $H_X(y) \leq k$. By the inductive hypothesis, we know that the survivor is caught after $time(y)$ rounds if the survivor remains in $G[\Lambda_X(y)]$. Therefore, the survivor must leave $G[\Lambda_X(y)]$ after $time(y) - 1 \leq d - 1$ steps, otherwise it is caught.

The survivor's walk thus looks like Figure 14. Each time it enters one of these subtrees, we assign the zombies $z_1, z_2, \ldots z_c$ to (specific) vertices in that subtree's component. Since $c$ is the maximum load of any child of $x$, we have a sufficient number of zombies. By the inductive hypothesis, this number of zombies suffices to either catch the survivor if the survivor remains in the component for $d$ steps or force the survivor out of the component of a child of $x$ and back into $C_x$ in at most $d - 1$ steps.

The survivor's walk length is therefore at most $(|C_x|(\delta(G)-1)+1) + (d-1)(|C_x|(\delta(G)-1)+1)$. The first term is the number of rounds the survivor can spend in $C_x$ until it is caught. For the second term, we note that the survivor can enter $G[\Lambda_X(y)]$ where $y$ is a child of $x$ at most $(|C_x|(\delta(G)-1)+1)$ times. Each time it enters $G[\Lambda_X(y)]$ it must return to a vertex in $C_x$ in $d-1$ rounds, otherwise the survivor is caught on the $d$th round. Therefore, we have that the survivor is caught after $(|C_x|(\delta(G)-1)+1) + (d-1)(|C_x|(\delta(G)-1)+1) \leq d((|C_x|(\delta(G)-1))+1) = time(x)$.

□

**Zombie strategy based on proof:** This inductive proof gives rise to the following zombie winning strategy. Let $X$ be a cut decomposition tree of $G$, with root $r$, such that $load(r) = load(G)$. Initially place *load*(*r*) zombies arbitrarily on vertices of *G*. Initially, all zombies are unassigned. Let the survivor's initial position be some vertex $v$ of *G*. Let $v \in C_x$ for some $x \in X$. Assign an arbitrary unassigned zombie to each vertex in $C_y$ for all $y \in \pi_X(x)$. When it is the zombies' turn to move, each unassigned zombie and each zombie that has reached its assigned vertex remains at its current location. An assigned zombie takes a step towards its assigned vertex if that vertex is on a shortest path from the zombie's current location to the survivor's current location. If on its turn, the survivor moves out of $C_x$ and to a vertex $u \in C_y$, then by construction of $X$, we have that $y$ must be either an ancestor or a descendant of $x$. If $y$ is a descendant of $x$, we assign an arbitrary unassigned zombie to each vertex in each $C_w$ where $w$ is a node that is on $\pi_X(y)$ and not on $\pi_X(x)$. Since we have *load*(*r*) zombies, there are a sufficient number of unassigned zombies available. On the other hand, if $y$ is an ancestor of $x$, then we unassign each zombie that was assigned to a vertex in each $C_w$ where $w$ is a node that is on $\pi_X(x)$ and not on $\pi_X(y)$. The survivor will be caught after *time*(*r*) + 1 rounds. The correctness of this strategy follows from the proof of Theorem 7.

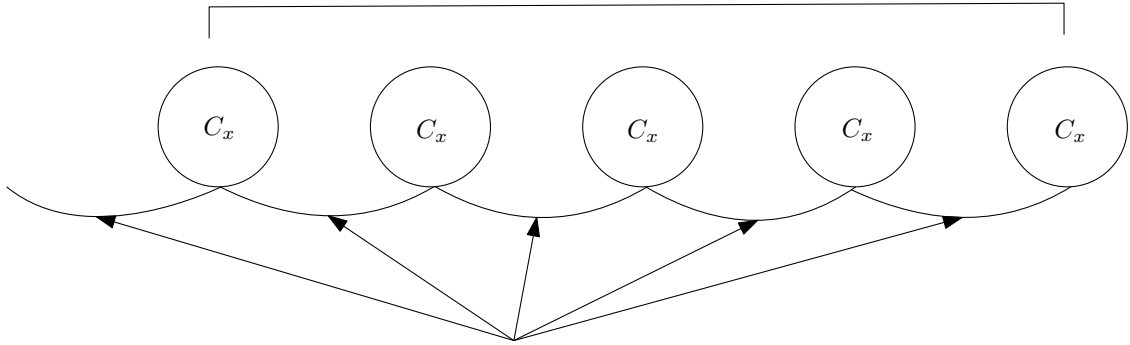

Figure 14: The path of the survivor

Among all cut decomposition trees realizing the load of *G*, we denote the value of the minimum height of such a tree by *cdh*(*G*). Among all cut decomposition trees whose load is *load*(*G*) and whose height is *cdh*(*G*), we denote the value of the minimum width among all such trees by *cdw*(*G*).

**Corollary 2.** *Given a connected graph G, $u_L(G) \leq load(G)$ and load(G) lazy zombies can catch the survivor in at most $(cdw(G)(\delta(G)-1)+1)^{cdh(G)+1}$ rounds.*

*Proof.* Let $X$ be a cut decomposition tree with root $r$, with load *load*(*G*), with height *cdh*(*G*) and with width *cdw*(*G*). Theorem 7 implies that $u_L(G) \leq load(G)$. Using Theorem 7 and the recursive definition of *time*, we show that *time*(*r*) is at most $(cdw(G)(\delta(G)-1)+1)^{cdh(G)+1}$.

Recall the definition of *time*(*x*) for a node $x \in X$. To simplify the notation, let $h = H_X$, and $a = cdw(G)(\delta(G)-1)+1$. We proceed by induction on the height of $X$ and prove that $time(r) \leq a^{h+1}$, for all $h \geq 0$. If $h$ is 0, then $time(r) = |C_r|(\delta(G)-1)+1 \leq a$. Assume that $time(r) \leq a^{h+1}$ when $h \leq k$, for $k \geq 0$. We now prove the inequality when $h = k+1$. By the definition of time, we have that $time(r) = \max_y time(y)(|C_r|(\delta(G)-1)+1) \leq \max_y time(y)a$, where $y$ is a child of $r$ in $X$. By the inductive hypothesis, we note that $\max_y time(y) \leq a^{h+1}$ since the height of $\Lambda_X(y)$ is at most $h$. This means that $time(r) \leq (a^{h+1})a \leq a^{h+2}$ as desired. Therefore, for all $h \geq 0$, we have $time(r) \leq a^{h+1}$. Substituting back, we get that $time(r)$ is $(cdw(G)(\delta(G)-1)+1)^{cdh(G)+1}$.

□

We will next show that for any graph $G$, $load(G)$ is exactly the treedepth $td(G)$. To do this, we will use the idea of a *compressed* version of a rooted tree. Given a rooted tree $T$, the compressed rooted tree $\tau(T)$ is a tree where every internal node has at least two children. To construct a compressed tree from $T$, we repeatedly apply the following process until there are no more internal nodes with exactly one child. Let $x_1, \ldots, x_i$ be a maximal path in $T$ where every node except $x_i$ has exactly one child, and $x_i$ has at least two children. Contract all the nodes in the path into a single node $s$ which we label as $\{x_1, \ldots, x_i\}$. The parent of $s$ in the tree is the parent of $x_1$ and the children of $s$ are the children of $x_i$.

**Theorem 8.** *For any connected graph G, $load(G) = td(G)$*

*Proof.* We begin by showing that $load(G) \leq td(G)$. Let $T$ be a rooted tree whose height is $td(G)$ such that $G \subseteq clos(T)$.

Recall that $V(T) = V(G)$. We construct a cut decomposition $(X, \mathcal{C})$ of $G$ from $T$. Let $X = \tau(T)$. Given a node $x \in X$, let the container $C_x$ be the union of the vertices of $T$ (i.e., vertices of $G$) that were contracted into $x$ during compression. In essence, the container is the label of the vertex as described in the compression process. We now show that the properties required for a valid cut-decomposition are satisfied:

- For every $v \in V(G)$, we need to show that there is a unique $x \in V(X)$ such that $v \in C_x$. Since $V(G) = V(T)$, $v$ can appear in only 1 container, since compressing paths does not place vertices in multiple containers.

- For every $uv \in E(G)$, we need to show that $\exists x, y \in V(X)$ such that $u \in C_x$, $v \in C_y$ and $x$ is an ancestor of $y$ in $X$. By the definition of $T$, we have that $u$ is an ancestor of $v$ in $T$. The property still holds in the compressed tree because if $u$ is an ancestor of $v$ prior to compression, then $x$ is an ancestor of $y$ after compression.

- For each non-leaf node $x \in X$, we need to show that $C_x$ is a cut-set of $G[\Lambda_X(x)]$. Since $X$ is a compressed tree, every internal vertex has at least 2 children. Let $y$ and $w$ be two such children of $x$. Removal of $x$ from $X$ decomposes the cut tree into at least 2 components, one containing $y$ and one containing $w$. The preceding property implies that $G[\Lambda_X(y)]$ and $G[\Lambda_X(w)]$ are two connected components when $C_x$ is removed from $G[\Lambda_X(x)]$. There cannot be an edge from a vertex in $G[\Lambda_X(y)]$ to a vertex in $G[\Lambda_X(w)]$ since $y$ is not an ancestor or descendant of $w$.

We have just shown that $(X, \mathcal{C})$ is a valid cut decomposition. The load of $(X, \mathcal{C})$ is the height of $T$ which implies that $load(G) \leq td(G)$.

We now show that $td(G) \leq load(G)$. Let $(X, \mathcal{C})$ be a cut decomposition of $G$. We will construct a rooted tree $T$ whose vertex set is $V(G)$ such that $G \subseteq clos(T)$.

To build $T$ from $X$, we essentially *uncompress* the nodes of $X$ as follows. For every node $x \in X$, we replace $x$ with a path towards the root consisting of the vertices in $C_x$. The first vertex in the path has an edge to the parent of $x$ and the last vertex in the path has edges to the children of $x$. Observe that this path induces a clique on the vertices of $C_x$ in $clos(T)$.

Consider an arbitrary edge $uv \in E(G)$. If $u$ and $v$ are in the same container $C_x$ of node $x$ in $X$, then $uv \in clos(T)$ since the vertices of $C_x$ are a path towards the root in $T$. If $u \in C_x$ and $v \in C_y$, then either $x$ is an ancestor of $y$ or $y$ is an ancestor of $x$ in $X$. This implies that $u$ is an ancestor of $v$ in $T$. Therefore, $uv \in clos(T)$. Thus, $G \subseteq clos(T)$.

Since the height of $T$ is the load of the tree decomposition, we conclude that $td(G) \leq load(G)$. □

For $0 < \alpha < 1$, a cut set $S \subseteq V(G)$ of $G$ is an $\alpha$-separator if every connected component of $G[V(G) - S]$ contains at most $\alpha n$ vertices. The size of the $\alpha$-separator is the cardinality of $S$. We highlight the relationship between the sizes of separators, treedepth, and treewidth below.

**Lemma 2.** *Let $G$ be a graph. Let $s_G : \{1, \ldots, n\} \to \mathbb{N}$ be a function defined as*

$$s_G(i) = \max_{A \subseteq V(G),\ |A| \leq i} \min\{|S| : S \text{ is a } \tfrac{1}{2}\text{-separator of } G[A]\}.$$

*The following inequalities hold:*

$$s_G(n) \leq load(G) = td(G) \leq \sum_{i=0}^{\log n} s_G(n/2^i) \leq (tw(G) + 1) \log n$$

*Proof.* The inequality $s_G(n) \leq td(G) \leq \sum_{i=0}^{\log n} s_G(n/2^i)$ is proven in Lemma 6.6 in [34]. Since it was shown by Robertson and Seymour [39] that $s_G(i) \leq tw(G) + 1$ for all $i \in [1, n]$, we have that $\sum_{i=0}^{\log n} s_G(n/2^i) \leq (tw(G) + 1) \log n$. The equality $load(G) = td(G)$ is proven in Theorem 8. □

$s_G(n)$ is sometimes called the separation number of $G$. The bound in Lemma 2 is tight in certain cases. For example, the treewidth of a path on $n$ vertices is 1 whereas the treedepth is $\Theta(\log n)$. However, for certain classes of graphs, we can remove the $\log n$ term on the upper bound in Lemma 2. Essentially, if $s_G(n/2^i) \leq cs_G(n)/2^i$ for some constant $c$, then $\sum_{i=0}^{\log n} s_G(n/2^i) \leq cs_G(n) \sum_{i=0}^{\log n} 1/2^i \leq 2cs_G(n) \leq 2c(tw(G) + 1)$. Thus, we have that $td(G)$ is $O(tw(G))$ in this case. Informally, this happens when the size of a separator for any subgraph of size $i$ is at most $i^\varepsilon$, for $0 < \varepsilon < 1$. This is summarized by the following:

**Corollary 3** (Corollary 6.2 in [34]). *Let $0 < \alpha < 1$, let $c > 0$ be a constant and let $\mathcal{G}$ be a hereditary class of graphs such that every $G \in \mathcal{G}$ with $n$ vertices has $tw(G) \leq cn^\alpha$, then every $G \in \mathcal{G}$ has $td(G) \leq \frac{c}{1-2^{-\alpha}} n^\alpha$.*

Treewidth, treedepth and separators are well-studied graph parameters. We highlight a few of the implications of our bound that $u_L(G) \leq td(G)$. The interested reader should consult the following comprehensive surveys on this topic [9, 34, 26, 19].

**Corollary 4.** *For connected planar graphs $G$, $u_L(G)$ is $O(\sqrt{n})$. These lazy zombies can catch the survivor in at most $n^{O(\log n)}$ rounds.*

*Proof.* The Planar Separator Theorem [31] implies that a planar graph $G$ has $tw(G)$ in $O(\sqrt{n})$. The bound on $load(G)$, and hence on $u_L(G)$, follows from Corollary 3 and Theorem 7.

An upper bound on the time of the cut decomposition is determined as follows. In Corollary 2, an upper bound of $O((cdw(G)(\delta(G)-1))^{cdh(G)+1})$ is proven for the time. For planar graphs, $cdw(G) \leq c\sqrt{n}$ for some constant $c$, $cdh(G) \leq \log_{3/2} n$, and $\delta(G) \leq n$. With these upper bounds, we conclude that the time is at most $n^{O(\log n)}$. □

With the genus-$g$ separator theorem of Gilbert, Hutchinson, and Tarjan [24], we can obtain:

**Corollary 5.** *There is a constant $c$ such that all connected genus-$g$ graphs $G$ have $u_L(G) \leq c\sqrt{gn}$. These lazy zombies can catch the survivor in at most $n^{O(\log n)}$ rounds.*

And similarly, with the excluded minor separation technique of Alon, Seymour, and Thomas [3], we get:

**Corollary 6.** *There is a constant $c$ such that all connected graphs $G$ with any excluded $h$-vertex minor $H$ have $u_L(G) \leq ch\sqrt{hn}$. These lazy zombies can catch the survivor in at most $n^{O(\log n)}$ rounds.*

Although $load(G)$ is an upper bound on $u_L(G)$, it is by no means tight. For example, if $G$ is a clique, then $load(G) = n$, whereas only one zombie suffices to catch a survivor in a clique. In an attempt to tighten some of these upper bounds, we try to leverage this idea that only one zombie suffices for a clique.

Given a cut decomposition $(X, \mathcal{C})$ of $G$, recall that our strategy is to assign a unique zombie to each vertex in a container. By assigning one zombie to each clique in a container of the cut decomposition rather than one zombie to each vertex in the container, we can improve the upper bound. For example, this immediately gives a bound of 1 for a clique which is tight. This idea leads to an alternative definition of load for a cut decomposition which we call $load^*$. In this definition, for $S \subseteq V(G)$, $\theta(S)$ is the clique cover number of the induced graph of $G$ on the vertices of $S$.

$$load^*(v) = \begin{cases} \theta(C_v) & \text{if } v \text{ is a leaf} \\ \theta(C_v) + \max_w load^*(w) & \text{otherwise, where } w \text{ is a child of } v \end{cases}$$

For some vertices and graphs, $\theta(C_v) = |C_v|$ (the cut set $C_v$ is an independent set), so without knowing more about the cuts, $load^*$ is no more useful than $load$. For some graphs, however, $load^*$ is a substantial improvement. The corresponding notion of time is:

$$time^*(v) = \begin{cases} \theta(C_v)\delta(G) + 1 & \text{if } v \text{ is a leaf} \\ (\theta(C_v)\delta(G) + 1)\max_w time^*(w) & \text{otherwise, where } w \text{ is a child of } v \end{cases}$$

**Theorem 9.** *Let $G$ be a connected graph with cut decomposition $(X, \mathcal{C})$. Suppose that $v$ is a vertex of $X$ and the survivor is restricted to the component of $v$. Then, $load^*(v)$ lazy zombies can capture the survivor in at most $time^*(v)$ rounds.*

*Proof.* The proof is analagous to the proof of Theorem 7 (including Lemma 1). The differences arise because we are now assigning zombies to cliques rather than vertices in $C_v$.

Suppose a zombie at location $z$ is assigned to a clique $K$ of a cut-decomposition container. Let the distance $d(z, K)$ from $z$ to $K$ be the minimum distance from $z$ to a vertex of $K$.

When the zombie is assigned to $K$, $d(z, K)$ can be at most $\delta(G)$. Each time the survivor is on a vertex of $K$, the zombie may decrease its distance to $K$: if the survivor is on a vertex $v$ of $K$ with $d(z, v) = d(z, K)$, then the zombie steps along a shortest path from $z$ to $v$. If the survivor is on a vertex $v$ of $K$ not at minimum distance to $z$, let $w$ be a vertex of $K$ such that $d(z, w)$ is minimum. Since $d(z, v) = d(z, w) + 1$, a step along a shortest path from $z$ to $w$ is also a step along a shortest path from $z$ to $v$ (and hence a legal zombie move). This decreases $d(z, K)$ by 1.

However, unlike in Theorem 7, where zombies attain their goal by becoming *adjacent* to their assigned vertex, the zombie here is in a position to block access to the entire clique only when $d(z, K) = 0$. This means that the survivor can make $\delta(G)$ moves to vertices of a clique before access to the clique is forbidden. □

**Corollary 7.** *Given a connected graph $G$, $u_L(G) \leq \text{load}^*(G) \leq \text{load}(G)$ and $\text{load}^*(G)$ lazy zombies can catch the survivor in at most $(\theta(G)(\delta(G)+1))^{cdh(G)+1}$ rounds.*

*Proof.* Similar to the proof of Corollary 2. □

## 6 Conclusion

We have shown that zombies are quite weak, in the worst case, and that a linear number of zombies are required to catch a survivor even in outerplanar graphs. We modified our lower bound construction to show that a linear number of zombies is still required in other related types of graphs such as the visibility graph of a simple polygon. We then showed that by simply allowing a zombie to be lazy and not have to move on its turn, only 2 lazy zombies are sufficient to catch a survivor in an outerplanar graph. Finally, we established that $k$ lazy zombies are sufficient to catch a survivor in graphs with treedepth $k$. We highlighted a few implications of this upper bound such as $(k+1)\log n$ lazy zombies are sufficient to catch a survivor on a graph with treewidth $k$.

Our linear lower bounds on zombies in outerplanar graphs is tight since a linear number of zombies is always sufficient to catch a survivor on any graph. However, our bound that $k$ lazy zombies is sufficient on graphs with treedepth $k$ is far from optimal since a clique has linear treedepth but 1 lazy zombie is sufficient. We attempted to address this issue by introducing the parameter $\text{load}^*(G)$. We showed that $\text{load}^*(G)$ lazy zombies are sufficient to catch a survivor in a graph of treedepth $k$ and that $\text{load}^*(G)$ is 1 when $G$ is a clique. However, this bound is still not tight and we leave as an open problem to find tighter upper bounds on the lazy zombie number of a graph.

### Acknowledgements

P. Bose would like to thank P. Morin for fruitful discussions on various aspects of structural graph theory and for bringing reference [26] to his attention.

## A  Appendix: Table and Figures for the Proof of Theorem 5

For reproducibility, the vertex locations of $Q$ (refer to Figure 6) are given in Table 2. The locations given are exact, not approximate.

In Figures 15 to 22), we show a few of the different steps in the survivor's initial walk through $Q_i$ (refer to the proof of Theorem 5). In each figure, it is the survivor's turn to play, the green vertex is the survivor's location, and the red vertices are the possible locations of zombies.

| | | | | | |
|---|---|---|---|---|---|
| s | (1.1263, 0.0172) | l5 | (-0.493, -0.4039) | m1 | (-0.1849, -0.538) |
| s3 | (1.1147, -0.0161) | j1 | (-0.5839, -0.3644) | m2 | (-0.1561, -0.5822) |
| s4 | (0.9363, -0.5145)) | j2 | (-0.7287, -0.0161) | m | (-0.1403, -0.6128) |
| b1 | (0.9097, -0.4961) | j | (-0.7403, 0.0172) | m3 | (-0.1391, -0.61) |
| b2 | (0.738, -0.0161) | j3 | (-0.752, -0.0161) | m4 | (-0.1288, -0.5918) |
| b | (0.7263, 0.0172) | j4 | (-0.8139, -0.1134) | e1 | (0.0536, -0.5973) |
| b3 | (0.713, -0.0161) | h1 | (-0.9259, -0.188) | e2 | (0.0787, -0.6005) |
| b4 | (0.5697, -0.2261) | h2 | (-1.127, -0.0161) | e | (0.1263, -0.6095) |
| d1 | (0.533, -0.2095) | h | (-1.1403, 0.0172) | e3 | (0.1397, -0.5828) |
| d2 | (0.3463, -0.0161) | g | (-1.207, -0.2495) | e4 | (0.2597, -0.3428) |
| d | (0.3263, 0.0172) | g3 | (-1.1823, -0.261) | c1 | (0.4063, -0.3561) |
| d3 | (0.3097, -0.0161) | g4 | (-1.1639, -0.2791) | c2 | (0.513, -0.5828) |
| d4 | (0.2697, -0.0728) | i1 | (-0.9602, -0.5739) | c | (0.5263, -0.6161) |
| f1 | (0.2197, -0.1095) | i2 | (-0.9496, -0.5901) | c3 | (0.5397, -0.5828) |
| f2 | (0.0677, -0.0532) | i | (-0.9403, -0.6161) | c4 | (0.6263, -0.4095) |
| f | (0.0597, -0.0495) | i3 | (-0.9362, -0.6091) | a1 | (0.8597, -0.4895) |
| f3 | (-0.0137, -0.2542) | i4 | (-0.8721, -0.5614) | a2 | (0.9097, -0.5828) |
| f4 | (-0.1356, -0.5384) | k1 | (-0.5895, -0.5365) | a | (0.9263, -0.6161) |
| l1 | (-0.1516, -0.5329) | k2 | (-0.557, -0.5828) | a3 | (0.943, -0.5828) |
| l2 | (-0.3363, 0.0051) | k | (-0.5403, -0.6161) | a4 | (0.963, -0.5695) |
| l | (-0.3403, 0.0172) | k3 | (-0.538, -0.6121) | t1 | (1.0063, -0.5595) |
| l3 | (-0.3513, -0.0176) | k4 | (-0.5035, -0.5652) | t2 | (1.0597, -0.5828) |
| l4 | (-0.4109, -0.2048) | k5 | (-0.355, -0.5358) | t | (1.1263, -0.6161) |

Table 2: Vertex locations of Q.

Figure 15: Round 1 in the survivor's walk.

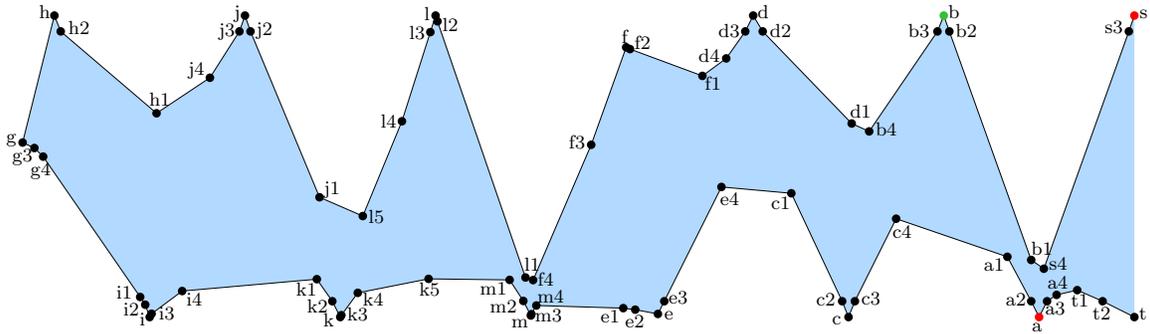

Figure 16: Round 2 in the survivor's walk.

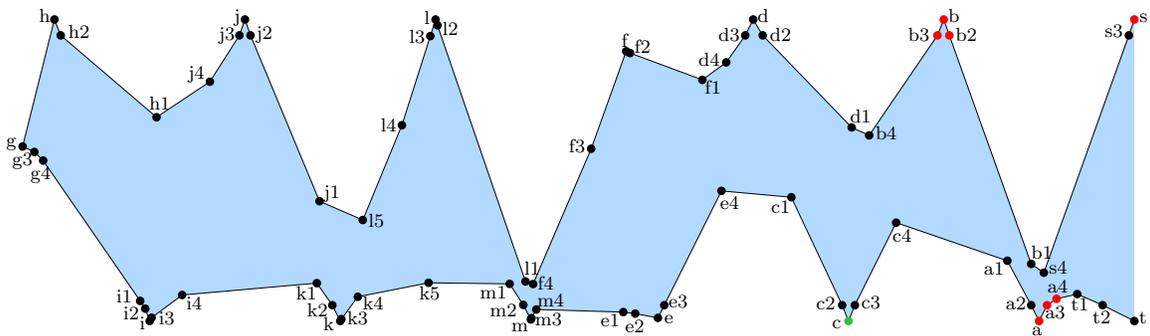

Figure 17: Round 3 in the survivor's walk.

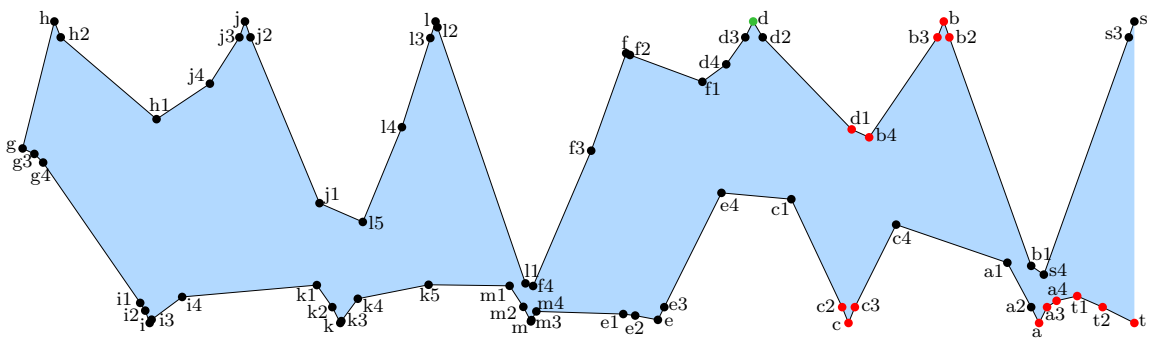

Figure 18: Round 4 in the survivor's walk.

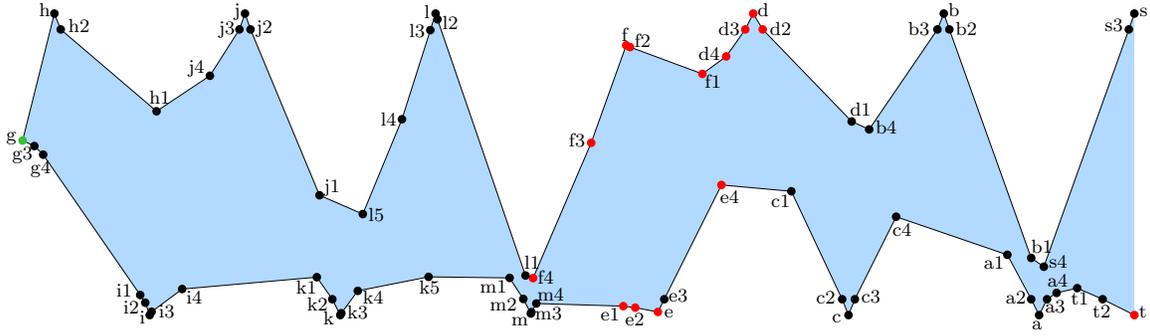

Figure 19: Round 6 in the survivor's walk.

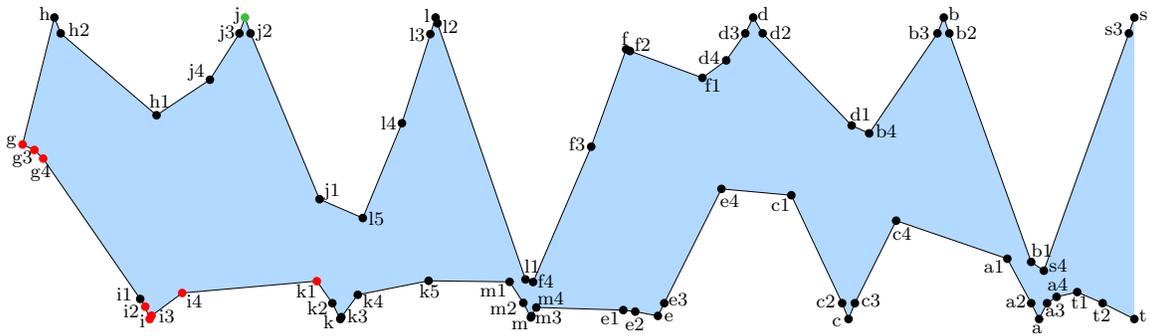

Figure 20: Round 9 in the survivor's walk.

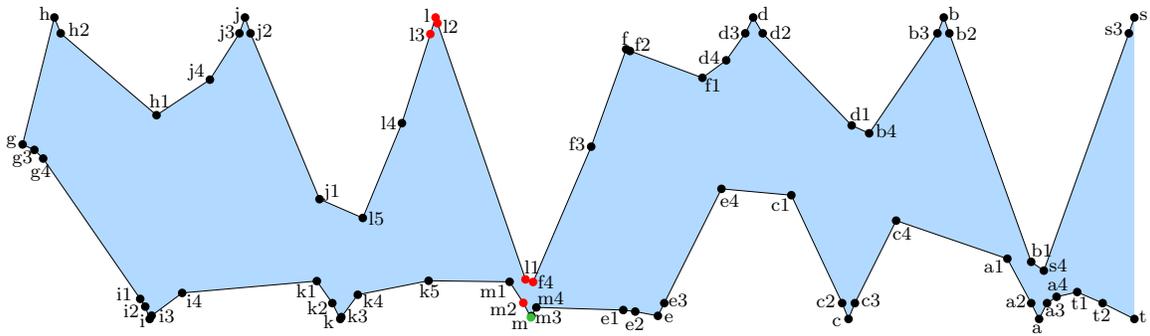

Figure 21: Round 12 in the survivor's walk.

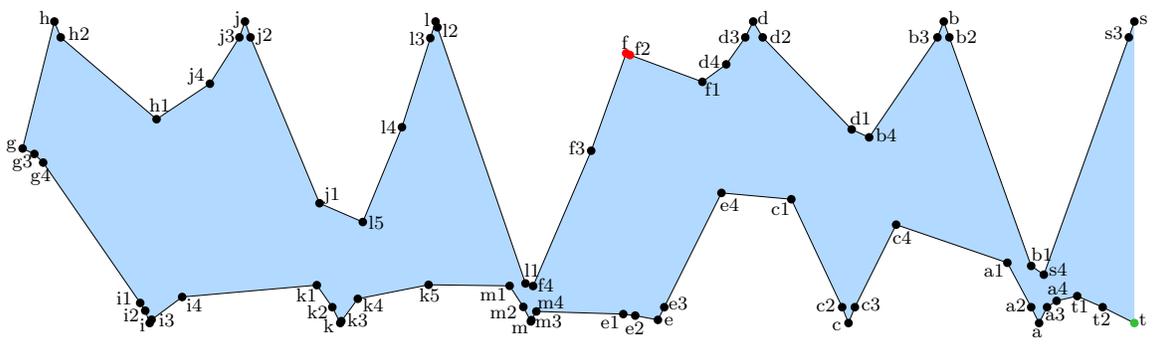

Figure 22: Round 14 in the survivor's walk.